\documentclass[12pt]{article}
\UseRawInputEncoding

\usepackage{amssymb,amsmath,amsthm,amsfonts,enumerate, bbm, mathdots}

\usepackage{latexsym}
\usepackage{amscd}
\usepackage{stmaryrd}
\usepackage{color}
\usepackage[all]{xypic}
\usepackage{epsfig}
\usepackage{graphics}
\usepackage{ifthen}
\usepackage{varioref}
\usepackage{rotating}
\usepackage{extarrows}
\usepackage{cite}
\usepackage{mathrsfs}

\numberwithin{equation}{section}

\theoremstyle{plain}
\newtheorem{thm}{Theorem}[section]
\newtheorem{cor}[thm]{Corollary}
\newtheorem{lem}[thm]{Lemma}
\newtheorem{prop}[thm]{Proposition}
\newtheorem{rem}[thm]{Remark}

\newmuskip\pFqmuskip
\newcommand*\pFq[6][8]{%
  \begingroup 
  \pFqmuskip=#1mu\relax
  \mathchardef\normalcomma=\mathcode`,
  \mathcode`\,=\string"8000
  \begingroup\lccode`\~=`\,
  \lowercase{\endgroup\let~}\pFqcomma
  {}_{#2}F_{#3}{\left(\genfrac..{0pt}{}{#4}{#5};#6\right)}%
  \endgroup
}
\newcommand{\pFqcomma}{{\normalcomma}\mskip\pFqmuskip}

\makeatletter
\newenvironment{proofof}[1]{\par
  \pushQED{\qed}%
  \normalfont \topsep6\p@\@plus6\p@\relax
  \trivlist
  \item[\hskip\labelsep
        \bfseries
    Proof of #1\@addpunct{.}]\ignorespaces
}{%
  \popQED\endtrivlist\@endpefalse
}
\makeatother

\definecolor{darkgreen}{rgb}{0.0625,0.64,0.0625}
\usepackage[
            pdfstartview=FitH,
            CJKbookmarks=true,
            bookmarksnumbered=true,
            bookmarksopen=true,
            colorlinks,
            pdfborder=001,
            linkcolor=blue,
            anchorcolor=green,
            citecolor=red
            ]{hyperref}

\newfont{\scyr}{wncyr10 scaled 550}

\def\proof{\noindent {\bf Proof.\;}}

\def\Ti{\operatorname{Ti}}
\allowdisplaybreaks

\begin{document}
\title{Relations of multiple $t$-values of general level\thanks{This work was supported by the Fundamental Research Funds for the Central Universities (grant number 22120210552).}}
\date{\small ~ \qquad\qquad School of Mathematical Sciences, Tongji University \newline No. 1239 Siping Road,
Shanghai 200092, China}

\author{Zhonghua Li\thanks{E-mail address: zhonghua\_li@tongji.edu.cn} ~and ~Zhenlu Wang\thanks{E-mail address: zlw@tongji.edu.cn}}
\maketitle
\begin{abstract}
We study the relations of multiple $t$-values of general level. The generating function of sums of multiple $t$-(star) values of level $N$ with fixed weight, depth and height is represented by the generalized hypergeometric function $_3F_2$, which generalizes the results for multiple zeta(-star) values and multiple $t$-(star) values. As applications, we obtain formulas for the generating functions of sums of multiple $t$-(star) values of level $N$ with height one and maximal height and a weighted sum formula for sums of multiple $t$-(star) values of level $N$ with fixed weight and depth. Using the stuffle algebra, we also get the symmetric sum formulas and Hoffman's restricted sum formulas for multiple $t$-(star) values of level $N$. Some evaluations of multiple $t$-star values of level $2$ with one-two-three indices are given.
\end{abstract}

{\small

}


\section{Introduction}

Let $\mathbf{k}=(k_1,\ldots,k_n)$ be a finite sequence of positive integers, which is called an index.  The quantities
$$k_1+\cdots+k_n,\quad  n \quad \text{and} \quad |\{j|1\leq j\leq n,k_j\geq 2\}|$$
are called the weight, depth and height of the index $\mathbf{k}$, respectively. The index $\mathbf{k}$ is called admissible if $k_1>1$.

The multiple zeta values of any level were introduced by H. Yuan and J. Zhao in \cite{Yuan-Zhao} and they focused on the double zeta values of level $N$, where $N$ is a fixed positive integer. Set $R=R_N=\mathbb{Z}/N\mathbb{Z}$. For an admissible index $\mathbf{k}=(k_1,\ldots,k_n)$ and any $\mathbf{a}=(a_1,a_2,\ldots,a_n)\in R^n$, the multiple zeta value of level $N$ is defined by
\begin{align*}
&\zeta_{N}(\mathbf{k};\mathbf{a})=\sum\limits_{m_1>\cdots>m_n>0\atop m_i\equiv a_i\pmod N}\frac{1}{m_1^{k_1}\cdots m_n^{k_n}}.
\end{align*}
Z. Li and Z. Wang studied the algebraic framework for the double shuffle relations of multiple zeta values of level $N$ and gave sum formulas and weighted sum formulas of double zeta values of level $2$ and $3$ in \cite{Li-Wang}. Two variants of multiple zeta values of level $2$ called the multiple $t$-values \cite{Hoffman19} and the multiple $T$-values \cite{Kaneko-Tsumura} have been studied recent years, see for example \cite{Charlton,Hoffman19,Kaneko-Tsumura,Li-Song,Murakami2021,Takeyama2021}.

In this paper, we consider multiple $t$-(star) values of level $N$, which generalize the multiple $t$-(star) values. For an admissible index $\mathbf{k}=(k_1,\ldots,k_n)$ and any $a\in \{1,2,\ldots,N\}$, we define the multiple $t$-value of level $N$ and the multiple $t$-star value of level $N$ respectively by
\begin{align*}
t_{N,a}(\mathbf{k})=t_{N,a}(k_1,\ldots,k_n)=\sum\limits_{m_1>\cdots>m_n>0\atop m_i\equiv a\pmod N}\frac{1}{m_1^{k_1}\cdots m_n^{k_n}},\\
t_{N,a}^{\star}(\mathbf{k})=t_{N,a}^{\star}(k_1,\ldots,k_n)=\sum\limits_{m_1\geq\cdots\geq m_n>0\atop m_i\equiv a\pmod N}\frac{1}{m_1^{k_1}\cdots m_n^{k_n}}.
\end{align*}
Note that if the index is empty, we treat the values $t_{N,a}(\varnothing)=t_{N,a}^\star(\varnothing)=1$. And we have the following iterated integral representations
\begin{align}
&t_{N,a}(\mathbf{k})=\int_0^1\left(\frac{dt}{t}\right)^{k_1-1}\frac{t^{N-1}dt}{1-t^N}\cdots\left(\frac{dt}{t}\right)^{k_{n-1}-1}\frac{t^{N-1}dt}{1-t^N}\left(\frac{dt}{t}\right)^{k_n-1}\frac{t^{a-1}dt}{1-t^N},\label{Eq:t_N,a integral}\\
&t_{N,a}^{\star}(\mathbf{k})=\int_0^1\left(\frac{dt}{t}\right)^{k_1-1}\frac{dt}{t(1-t^N)}\cdots\left(\frac{dt}{t}\right)^{k_{n-1}-1}\frac{dt}{t(1-t^N)}\left(\frac{dt}{t}\right)^{k_n-1}\frac{t^{a-1}dt}{1-t^N},\label{Eq:t_N,a^star integral}
\end{align}
where for one forms $w_i=f_i(t)dt$, $i=1,2,\ldots,k$, we define that
\begin{align*}
\int_0^zw_1w_2\cdots w_k=\int\limits_{z>t_1>t_2>\ldots>t_k>0}f_1(t)f_2(t)\cdots f_k(t)dt_1dt_2\cdots dt_k.
\end{align*}
Obviously,
\begin{align*}
&t_{1,1}(\mathbf{k})=\zeta(\mathbf{k})=\sum\limits_{m_1>\cdots>m_n>0}\frac{1}{m_1^{k_1}\cdots m_n^{k_n}},\\
&t_{1,1}^{\star}(\mathbf{k})=\zeta^{\star}(\mathbf{k})=\sum\limits_{m_1\geq\cdots\geq m_n>0}\frac{1}{m_1^{k_1}\cdots m_n^{k_n}},
\end{align*}
which are the multiple zeta value and the multiple zeta-star value \cite{Hoffman92} respectively. We also have $t_{N,N}(\mathbf{k})=N^{-k}\zeta(\mathbf{k})$ and $t_{N,N}^\star(\mathbf{k})=N^{-k}\zeta^{\star}(\mathbf{k})$, where $k$ is the weight of the index $\mathbf{k}$. Taking $N=2$ and $a=1$, we get the multiple $t$-value
$$t_{2,1}(\mathbf{k})=t(\mathbf{k})=\sum\limits_{m_1>\cdots>m_n>0\atop m_i:\text{odd}}\frac{1}{m_1^{k_1}\cdots m_n^{k_n}},$$
and the multiple $t$-star value
$$t_{2,1}^{\star}(\mathbf{k})=t^{\star}(\mathbf{k})=\sum\limits_{m_1\geq\cdots\geq m_n>0\atop m_i:\text{odd}}\frac{1}{m_1^{k_1}\cdots m_n^{k_n}}.$$

Y. Ohno and D. Zagier \cite{Ohno-Zagier} studied the generating function of sums of multiple zeta values with fixed weight, depth and height by using the Gaussian hypergeometric function. T. Aoki, Y. Kombu and Y. Ohno \cite{Aoki-Kombu-Ohno} represented the generating function of sums of multiple zeta-star values with fixed weight, depth and height by the generalized hypergeometric function $_3F_2$. Here for a positive integer $m$ and complex numbers $b_1,\ldots,b_{m+1}$, $c_1,\ldots,c_m$ with none of $c_i$ is zero or a negative integer, the generalized hypergeometric function $_{m+1}F_m$ is defined by
\begin{align*}
\pFq{m+1}{m}{b_1,\ldots,b_{m+1}}{c_1,\ldots,c_m}{z}=\sum\limits_{n=0}^\infty\frac{(b_1)_n\cdots(b_{m+1})_n}{(c_1)_n\cdots(c_m)_n}\frac{z^n}{n!},
\end{align*}
where the Pochhammer symbol $(b)_n$ is defined by
$$(b)_n=\frac{\Gamma(b+n)}{\Gamma(b)}=\begin{cases}
1 & \text{if}\quad n=0, \\
b(b+1)\cdots(b+n-1) & \text{if}\quad n>0.
\end{cases}$$
The above series is absolutely and uniformly convergent for $|z|<1$ and the convergence also extends over the unit circle if $\Re\left(\sum c_i-\sum b_i\right)>0$. As analogues of Ohno-Zagier relations, the Ohno-Zagier type relations for the multiple $T$-values were studied by Y. Takeyama in \cite{Takeyama2021} and the Ohno-Zagier type relations for the multiple $t$-(star) values were given by Z. Li and Y. Song in \cite{Li-Song}. Similarly as in \cite{Li-Song}, we study the Ohno-Zagier type relations for multiple $t$-(star) values of level $N$ and give some applications in Section \ref{Sec:Ohno-Zagier type relation}.

Due to the infinite series representations, it is easy to see that the multiple $t$-(star) values of level $N$ satisfy stuffle relations. Hence we obtain the symmetric sum formulas and Hoffman's restricted sum formulas for multiple $t$-(star) values of level $N$ by establishing the algebraic setup in Section \ref{Sec:Sum formulas}.

D. Zagier found evaluations of $\zeta(2,\ldots,2,3,2,\ldots,2)$ and $\zeta^\star(2,\ldots,2,3,2,\ldots,2)$ by establishing the generating functions in \cite{Zagier2012}. By similar method, T. Murakami  gave the evaluation of $t(2,\ldots,2,3,2,\ldots,2)$ in \cite{Murakami2021} and S. Charlton  provided the evaluation of $t(2,\ldots,2,1,2,\ldots,2)$ in \cite{Charlton}. We get evaluations of their star version $t^\star(2,\ldots,2,3,2,\ldots,2)$ and $t^\star(2,\ldots,2,1,2,\ldots,2)$ in Section \ref{Sec:Evaluations}.

\section{Ohno-Zagier type relations}\label{Sec:Ohno-Zagier type relation}

For nonnegative integers $k,n,s$, let $I_0(k,n,s)$ be the set of admissible indices of weight $k$, depth $n$ and height $s$. Notice that $I_0(k,n,s)$ is nonempty if and only if $k\geq n+s$ and $n\geq s\geq1$. Define the sums
\begin{align*}
&G_{0,N,a}(k,n,s)=\sum\limits_{\mathbf{k}\in I_0(k,n,s)}t_{N,a}(\mathbf{k}),\\
&G_{0,N,a}^{\star}(k,n,s)=\sum\limits_{\mathbf{k}\in I_0(k,n,s)}t_{N,a}^{\star}(\mathbf{k}).
\end{align*}
Then we obtain Ohno-Zagier type relations which represent the generating functions of $G_{0,N,a}(k,n,s)$ and $G_{0,N,a}^{\star}(k,n,s)$ by the generalized hypergeometric function $_3F_2$.

\begin{thm}\label{Thm:t_N}
For formal variables $u,v,w$,
\begin{align*}
\sum\limits_{k\geq n+s,n\geq s\geq1}G_{0,N,a}(k,n,s)u^{k-n-s}v^{n-s}w^{s-1}=\frac{1}{a(a-u)}\pFq{3}{2}{\alpha,\beta,1}{\frac{a+N}{N},\frac{a+N-u}{N}}{1},
\end{align*}
where $\alpha,\beta$ are determined by $\alpha+\beta=\frac{1}{N}(2a-u+v)$ and $\alpha\beta=\frac{1}{N^2}(a(a-u+v)-uv+w)$.
\end{thm}

\begin{thm}\label{Thm:t-star_N}
For formal variables $u,v,w$,
\begin{align*}
\sum\limits_{k\geq n+s,n\geq s\geq1}G_{0,N,a}^{\star}(k,n,s)u^{k-n-s}v^{n-s}w^{s-1}=\frac{1}{a(a-u-v)+uv-w}\pFq{3}{2}{\frac{a}{N},\frac{a-u}{N},1}{\alpha^{\star},\beta^{\star}}{1},
\end{align*}
where $\alpha^{\star},\beta^{\star}$ are determined by $\alpha^{\star}+\beta^{\star}=\frac{1}{N}(2a+2N-u-v)$ and $\alpha^{\star}\beta^{\star}=\frac{1}{N^2}((a+N)(a+N-u-v)+uv-w)$.
\end{thm}

From the above theorems and using some summation formulas of $_3F_2$, we obtain some corollaries in Subsection \ref{Subsec:Applications}. For example, Theorem \ref{Thm:t_N} and Theorem \ref{Thm:t-star_N} deduce the Ohno-Zagier type relations for multiple zeta values and multiple zeta-star values respectively. Furthermore, we study the generating functions of sums of multiple $t$-(star) values of level $N$ with height one and maximal height and obtain a weighted sum formula of sums of multiple $t$-(star) values of level $N$ with fixed weight, fixed depth and $N=2a$. We will prove the above theorems in Subsection \ref{Subsec:Proofs}.

\subsection{Applications}\label{Subsec:Applications}

\subsubsection{Level one and level two cases}

Setting $N=a=1$ or only $a=N$ in Theorem \ref{Thm:t_N} and Theorem \ref{Thm:t-star_N}, we can obtain the Ohno-Zagier type relations for the multiple zeta(-star) values.

Setting $N=a=1$ in Theorem \ref{Thm:t_N}, we get $\alpha+\beta=2-u+v$ and $\alpha\beta=1-u+v-uv+w$. Let $\alpha_1=u+\alpha-1$ and $\beta_1=u+\beta-1$, then $\alpha_1+\beta_1=u+v$ and $\alpha_1\beta_1=w$. Hence we have
\begin{align*}
&\frac{1}{1-u}\pFq{3}{2}{\alpha,\beta,1}{2,2-u}{1}=\frac{1}{uv-w}\left[1-\pFq{2}{1}{\alpha-1,\beta-1}{1-u}{1}\right]\\
&=\frac{1}{uv-w}\left[1-\pFq{2}{1}{\alpha_1-u,\beta_1-u}{1-u}{1}\right]\\
&=\frac{1}{uv-w}\left[1-\frac{\Gamma\left(1-u\right)\Gamma\left(1-v\right)}{\Gamma(1-\alpha_1)\Gamma(1-\beta_1)}\right].
\end{align*}
Using the expansion
\begin{align}\label{Eq:Gamma's expansion}
\Gamma(1-z)=\exp\left(\gamma z+\sum\limits_{n=2}^\infty \frac{\zeta(n)}{n}z^n\right),
\end{align}
we obtain the following result, which is the Ohno-Zagier relation for the multiple zeta values.

\begin{cor}[{\cite[Theorem 1]{Ohno-Zagier}}]
For formal variables $u,v,w$,
\begin{align*}
&\sum\limits_{k\geq n+s,n\geq s\geq1}\left(\sum\limits_{\mathbf{k}\in I_0(k,n,s)}\zeta(\mathbf{k})\right)u^{k-n-s}v^{n-s}w^{s-1}\\
&=\frac{1}{uv-w}\left\{1-\exp\left(\sum\limits_{n=2}^\infty \frac{\zeta(n)}{n}(u^n+v^n-\alpha^n-\beta^n)\right)\right\}.
\end{align*}
where $\alpha$ and $\beta$ are determined by $\alpha+\beta=u+v$ and $\alpha\beta=w$.
\end{cor}

Similarly, setting $N=a=1$ in Theorem \ref{Thm:t-star_N}, we get $\alpha^\star+\beta^\star=4-u-v$ and $\alpha^\star\beta^\star=4-2(u+v)+uv-w$. Let $\alpha_1^\star=2-\alpha^\star$ and $\beta_1^\star=2-\beta^\star$, then $\alpha_1^\star+\beta_1^\star=u+v$ and $\alpha_1^\star\beta_1^\star=uv-w$. Using the summation formula \cite[7.4.4.1]{Prudnikov-Brychkov-Marichev}
\begin{align*}
\pFq{3}{2}{a,b,c}{d,e}{1}=\frac{\Gamma(d)\Gamma(d+e-a-b-c)}{\Gamma(d-c)\Gamma(d+e-a-b)}\pFq{3}{2}{e-a,e-b,c}{e,d+e-a-b}{1}
\end{align*}
with $a=c=1$, $b=1-u$, $d=\alpha^\star$ and $e=\beta^\star$, we have
\begin{align*}
&\frac{1}{1-u-v+uv-w}\pFq{3}{2}{1,1,1-u}{\alpha^\star,\beta^\star}{1}\\&=\frac{1}{1-u-v+uv-w}\frac{\alpha^\star-1}{1-v}\pFq{3}{2}{\beta^\star-1,\beta^\star-1+u,1}{\beta^\star,2-v}{1}\\
&=\frac{1}{(1-\beta_1^\star)(1-v)}\pFq{3}{2}{1-\beta_1^\star,1-\beta_1^\star+u,1}{2-\beta_1^\star,2-v}{1}.
\end{align*}
Then we obtain the following result, which is the Ohno-Zagier type relation for the multiple zeta-star values.

\begin{cor}[{\cite[Proposition 3.1]{Aoki-Kombu-Ohno}}]
For formal variables $u,v,w$,
\begin{align*}
&\sum\limits_{k\geq n+s,n\geq s\geq1}\left(\sum\limits_{\mathbf{k}\in I_0(k,n,s)}\zeta^\star(\mathbf{k})\right)u^{k-n-s}v^{n-s}w^{s-1}\\&=\frac{1}{(1-\beta^\star)(1-v)}\pFq{3}{2}{1-\beta^\star,1-\beta^\star+u,1}{2-\beta^\star,2-v}{1},
\end{align*}
where $\alpha^\star$ and $\beta^\star$ are determined by $\alpha^\star+\beta^\star=u+v$ and $\alpha^\star\beta^\star=uv-w$.
\end{cor}

For the level two case, setting $N=2$ and $a=1$ in Theorem \ref{Thm:t_N} and Theorem \ref{Thm:t-star_N}, we get the Ohno-Zagier relations for multiple $t$-(star) values \cite[Theorems 1.1 and 1.2]{Li-Song}.


\subsubsection{Sum of height one}

To save space, denote by $\{k\}^n$ the sequence of $k$ repeated $n$ times.

Setting $w=0$ in Theorem \ref{Thm:t_N}, we have $\alpha=\frac{a-u}{N}$, $\beta=\frac{a+v}{N}$. Therefore, we get the generating function of height one multiple $t$-values of level $N$.

\begin{cor}
For formal variables $u,v$, we have
\begin{align*}
\sum\limits_{k\geq n+1,n\geq1}t_{N,a}(k-n+1,\{1\}^{n-1})u^{k-n-1}v^{n-1}=\frac{1}{a(a-u)}\pFq{3}{2}{\frac{a-u}{N},\frac{a+v}{N},1}{\frac{a+N}{N},\frac{a+N-u}{N}}{1}.
\end{align*}
\end{cor}

Setting $w=0$ in Theorem \ref{Thm:t-star_N}, we have $\alpha^{\star}=\frac{a+N-u}{N}$, $\beta^{\star}=\frac{a+N-v}{N}$. Hence, we get the generating function of heignt one multiple $t$-star values of level $N$.

\begin{cor}
For formal variables $u,v$, we have
\begin{align*}
\sum\limits_{k\geq n+1,n\geq1}t_{N,a}^{\star}(k-n+1,\{1\}^{n-1})u^{k-n-1}v^{n-1}=\frac{1}{(a-u)(a-v)}\pFq{3}{2}{\frac{a}{N},\frac{a-u}{N},1}{\frac{a+N-u}{N},\frac{a+N-v}{N}}{1}.
\end{align*}
\end{cor}

Similar to \cite[Lemma 5.2]{Hoffman19}, we have the following lemma.

\begin{lem}\label{Lem:3F2}
For integers $a,N$ with $1\leq a\leq N$ and a formal variable $x$, if
\begin{align*}
\frac{x}{a-x}\pFq{3}{2}{p,q,\frac{a-x}{N}}{s,\frac{a+N-x}{N}}{1}=\sum\limits_{j=1}^\infty t_jx^j,
\end{align*}
then
\begin{align*}
t_n=a^{-n}\pFq{n+2}{n+1}{p,q,\{\frac{a}{N}\}^n}{s,\{\frac{a+N}{N}\}^n}{1}.
\end{align*}
\end{lem}

\proof
By definition,
\begin{align}\label{Eq:3F2}
\frac{x}{a-x}\pFq{3}{2}{p,q,\frac{a-x}{N}}{s,\frac{a+N-x}{N}}{1}=x\sum\limits_{i=0}^\infty \frac{(p)_i(q)_i}{(s)_ii!}\frac{1}{a+Ni-x}.
\end{align}
Since the $n$th derivative of the right-hand side of \eqref{Eq:3F2} with respect to $x$ is
\begin{align*}
n!\sum\limits_{i=0}^\infty \frac{(p)_i(q)_i}{(s)_ii!}\frac{1}{(a+Ni-x)^n}+n!x\sum\limits_{i=0}^\infty \frac{(p)_i(q)_i}{(s)_ii!}\frac{1}{(a+Ni-x)^{n+1}},
\end{align*}
we have
\begin{align*}
t_n=\frac{1}{n!}\left.\dfrac{d^n}{dx^n}\right|_{x=0}\frac{x}{a-x}\pFq{3}{2}{p,q,\frac{a-x}{N}}{s,\frac{a+N-x}{N}}{1}=\sum\limits_{i=0}^\infty \frac{(p)_i(q)_i}{(s)_ii!}\frac{1}{(a+Ni)^n}.
\end{align*}
Therefore, the conclusion follows since $\left(\frac{a}{N}\right)_i/\left(\frac{a+N}{N}\right)_i=\frac{a}{a+Ni}$.
\qed

Using Lemma \ref{Lem:3F2}, for any integer $m\geq2$, we have
\begin{align*}
\sum\limits_{n=1}^\infty t_{N,a}(m,\{1\}^{n-1})v^{n-1}&=a^{-m}\pFq{m+1}{m}{1,\frac{a+v}{N},\{\frac{a}{N}\}^{m-1}}{\{\frac{a+N}{N}\}^m}{1},\\
\sum\limits_{n=1}^\infty t_{N,a}^{\star}(m,\{1\}^{n-1})v^{n-1}&=\frac{1}{a^{m-1}(a-v)}\pFq{m+1}{m}{1,\{\frac{a}{N}\}^{m}}{\frac{a+N-v}{N},\{\frac{a+N}{N}\}^{m-1}}{1}.
\end{align*}


\subsubsection{Sum of maximal height}

We recall the following expansion formula of the gamma function with a proof for the convenience of the readers.

\begin{lem}\label{Lem:Gamma representation}
For integers $a,N$ with $1\leq a\leq N$, we have
\begin{align}
\Gamma\left(\frac{a-z}{N}\right)=\Gamma\left(\frac{a}{N}\right)\exp\left\{-\frac{1}{N}\psi\left(\frac{a}{N}\right)z+\sum\limits_{n=2}^\infty \frac{t_{N,a}(n)}{n}z^n\right\},
\end{align}
where $\psi$ denotes the digamma function.
\end{lem}

\proof
Since $\log\Gamma(z)=-\gamma z-\log z+\sum\limits_{k=1}^\infty \left(\frac{z}{k}-\log\left(1+\frac{z}{k}\right)\right)$ and $\log\Gamma(z+1)=\log\Gamma(z)+\log z$, we have
\begin{align*}
\frac{d}{dz}\log\Gamma(z+1)=-\gamma+\sum\limits_{k=1}^\infty\frac{z}{k(z+k)}
\end{align*}
and
\begin{align*}
\frac{d^2}{dz^2}\log\Gamma(z+1)=\sum\limits_{k=1}^\infty\frac{1}{(z+k)^2}.
\end{align*}
Therefore, we obtain that
\begin{align*}
\frac{d^2}{dz^2}\log\Gamma\left(\frac{a-N-z}{N}+1\right)&=\frac{1}{N^2}\left.\dfrac{d^2}{d\omega^2}\right|_{\omega=\frac{a-N-z}{N}}\log\Gamma(\omega+1)\\
&=\sum\limits_{n=2}^\infty(n-1)t_{N,a}(n)z^{n-2}.
\end{align*}
Hence, we have
\begin{align*}
\Gamma\left(\frac{a-z}{N}\right)=\exp\left\{c_2+c_1z+\sum\limits_{n=2}^\infty \frac{t_{N,a}(n)}{n}z^n\right\},
\end{align*}
where $c_2=\log\Gamma(\frac{a}{N})$ and
\begin{align*}
c_1=\frac{\gamma}{N}+\frac{1}{N}\sum\limits_{k=1}^\infty\frac{N-a}{k((k-1)N+a)}.
\end{align*}
Note that
\begin{align*}
\frac{1}{N}\sum\limits_{k=1}^m\frac{N-a}{k((k-1)N+a)}&=\sum\limits_{k=1}^m\left[\frac{1}{(k-1)N+a}-\frac{1}{kN}\right]\\
&=\frac{1}{N}\left[\psi\left(m+\frac{a}{N}\right)-\psi\left(\frac{a}{N}\right)-H_m\right],
\end{align*}
where $H_m$ denotes the $m$-th harmonic number. Since $H_m=\psi(m+1)+\gamma$, using the asymptotics of the digamma function, we find that
\begin{align*}
\frac{1}{N}\sum\limits_{k=1}^\infty\frac{N-a}{k((k-1)N+a)}=-\frac{\gamma}{N}-\frac{1}{N}\psi\left(\frac{a}{N}\right).
\end{align*}
Thus, we prove the result.
\qed

\begin{cor}
For formal variables $u,w$, we have
\begin{align}\label{Eq:G_0_maximal height representation}
1+\sum\limits_{k\geq 2n,n\geq1}G_{0,N,a}(k,n,n)u^{k-2n}w^n=\exp\left\{\sum\limits_{n=2}^\infty\frac{t_{N,a}(n)}{n}(u^n-x^n-y^n)\right\},
\end{align}
where $x,y$ are determined by $x+y=u$ and $xy=w$.
\end{cor}

\proof
Setting $v=0$ in Theorem \ref{Thm:t_N}, we get $\alpha+\beta=\frac{1}{N}(2a-u)$ and $\alpha\beta=\frac{1}{N^2}(a(a-u)+w)$. Let $x=a-N\alpha$ and $y=a-N\beta$, then $x+y=u$ and $xy=w$. Using the summation formula \cite[7.4.4.28]{Prudnikov-Brychkov-Marichev}
\begin{align}\label{Eq:3F2 maximal height}
\pFq{3}{2}{a,b,1}{c,2+a+b-c}{1}=\frac{1+a+b-c}{(1+a-c)(1+b-c)}\left(1-c+\frac{\Gamma(c)\Gamma(1+a+b-c)}{\Gamma(a)\Gamma(b)}\right)
\end{align}
with $a=\alpha$, $b=\beta$ and $c=\frac{a+N}{N}$, we find that
\begin{align*}
\frac{1}{a(a-u)}\pFq{3}{2}{\alpha,\beta,1}{\frac{a+N}{N},\frac{a+N-u}{N}}{1}=-\frac{1}{w}+\frac{1}{w}\frac{\Gamma\left(\frac{a}{N}\right)\Gamma\left(\frac{a-u}{N}\right)}{\Gamma\left(\frac{a-x}{N}\right)\Gamma\left(\frac{a-y}{N}\right)}.
\end{align*}
Then it is easy to finish the proof by Lemma \ref{Lem:Gamma representation}.
\qed

\begin{cor}
For formal variables $u,w$, we have
\begin{align}\label{Eq:G_0^star_maximal height representation}
1+\sum\limits_{k\geq 2n,n\geq1}G_{0,N,a}^\star(k,n,n)u^{k-2n}w^n=\exp\left\{\sum\limits_{n=2}^\infty\frac{t_{N,a}(n)}{n}((x^\star)^n+(y^\star)^n-u^n)\right\},
\end{align}
where $x^\star,y^\star$ are determined by $x^\star+y^\star=u$ and $x^\star y^\star=-w$.
\end{cor}
\proof
Setting $v=0$ in Theorem \ref{Thm:t-star_N}, we get $\alpha^\star+\beta^\star=\frac{1}{N}(2a+2N-u)$ and $\alpha^\star\beta^\star=\frac{1}{N^2}((a+N)(a+N-u)-w)$. Let $x^\star=a+N-N\alpha^\star$ and $y^\star=a+N-N\beta^\star$, then $x^\star+y^\star=u$ and $x^\star y^\star=-w$. Using \eqref{Eq:3F2 maximal height} with $a=\frac{a}{N}$, $b=\frac{a-u}{N}$ and $c=\alpha^\star$, we find that
\begin{align*}
\frac{1}{a(a-u)-w}\pFq{3}{2}{\frac{a}{N},\frac{a-u}{N}}{\alpha^\star,\beta^\star}{1}=-\frac{1}{w}+\frac{1}{w}\frac{\Gamma\left(\frac{a-x}{N}\right)\Gamma\left(\frac{a-y}{N}\right)}{\Gamma\left(\frac{a}{N}\right)\Gamma\left(\frac{a-u}{N}\right)},
\end{align*}
which with Lemma \ref{Lem:Gamma representation} implies the result.
\qed

By \eqref{Eq:G_0_maximal height representation} and \eqref{Eq:G_0^star_maximal height representation}, we find that
\begin{align*}
&\left(1+\sum\limits_{k\geq 2n,n\geq1}G_{0,N,a}(k,n,n)u^{k-2n}w^n\right)\\
&\quad\times \left(1+\sum\limits_{k\geq 2n,n\geq1}(-1)^nG_{0,N,a}^\star(k,n,n)u^{k-2n}w^n\right)=1.
\end{align*}
Furthermore, setting $u=0$ in \eqref{Eq:G_0_maximal height representation} and \eqref{Eq:G_0^star_maximal height representation}, we have
\begin{align*}
1+\sum\limits_{n=1}^\infty t_{N,a}(\{2\}^n)w^n&=\exp\left(\sum\limits_{n=1}^\infty\frac{(-1)^{n-1}t_{N,a}(2n)}{n}w^n\right),\\
1+\sum\limits_{n=1}^\infty t_{N,a}^\star(\{2\}^n)w^n&=\exp\left(\sum\limits_{n=1}^\infty\frac{t_{N,a}(2n)}{n}w^n\right).
\end{align*}
Using the symmetric sum formulas showed in Proposition \ref{Prop:SymSum}, we will give more general formulas in Corollary \ref{Cor:Rep-Arg}.


\subsubsection{A weighted sum formula}

Let $I_0(k,n)$ be the set of admissible indices of weight $k$ and depth $n$. Define
\begin{align*}
G_{0,N,a}(k,n)=\sum\limits_{\mathbf{k}\in I_0(k,n)}t_{N,a}(\mathbf{k}),\quad
G_{0,N,a}^\star(k,n)=\sum\limits_{\mathbf{k}\in I_0(k,n)}t_{N,a}^\star(\mathbf{k}).
\end{align*}

\begin{prop}
For a formal variable $u$, we have
\begin{align}\label{Eq:G_0,2a,a}
&\sum\limits_{k=2}^\infty\left(\sum\limits_{n=1}^{k-1}2^{n-1}G_{0,2a,a}(k,n)\right)u^{k-2}=\sum\limits_{k=2}^\infty\left(\sum\limits_{n=1}^{k-1}(-1)^{k-n-1}2^{n-1}G_{0,2a,a}^\star(k,n)\right)u^{k-2}\notag\\
&=a^{-2}2^{\frac{2u}{a}}t(2)\exp\left(\sum\limits_{n=2}^\infty\frac{4(2^{n-1}-1)}{na^n(2^n-1)}t(n)u^n\right).
\end{align}
\end{prop}

\proof
Let $uv=w$ in Theorem \ref{Thm:t_N}, we have $\alpha=\frac{a-u+v}{N}$, $\beta=\frac{a}{N}$. Hence, we get the generating function of sums of multiple $t$-values of level $N$ with fixed weight and depth:
\begin{align*}
\sum\limits_{k>n\geq1}G_{0,N,a}(k,n)u^{k-n-1}v^{n-1}=\frac{1}{a(a-u)}\pFq{3}{2}{\frac{a-u+v}{N},\frac{a}{N},1}{\frac{a+N-u}{N},\frac{a+N}{N}}{1}.
\end{align*}
Setting $v=2u$ and $N=2a$, we have
\begin{align*}
\sum\limits_{k>n\geq1}2^{n-1}G_{0,2a,a}(k,n)u^{k-2}=\frac{1}{a(a-u)}\pFq{3}{2}{\frac{a+u}{2a},\frac{1}{2},1}{\frac{3a-u}{2a},\frac{3}{2}}{1}.
\end{align*}
Using the summation formula \cite[7.4.4.21]{Prudnikov-Brychkov-Marichev}
\begin{align*}
&\pFq{3}{2}{a,b,c}{1+a-b,1+a-c}{1}\\
&=\frac{\sqrt{\pi}}{2^a}\frac{\Gamma(1+a-b)\Gamma(1+a-c)\Gamma\left(1+\frac{a}{2}-b-c\right)}{\Gamma\left(\frac{1+a}{2}\right)\Gamma\left(1+\frac{a}{2}-b\right)\Gamma\left(1+\frac{a}{2}-c\right)\Gamma(1+a-b-c)},\quad\Re(a-2b-2c)>-2,
\end{align*}
with $a=1$, $b=\frac{1}{2}$ and $c=\frac{a+u}{2a}$, we get
\begin{align*}
\frac{1}{a(a-u)}\pFq{3}{2}{\frac{a+u}{2a},\frac{1}{2},1}{\frac{3a-u}{2a},\frac{3}{2}}{1}&=\frac{1}{a(a-u)}\frac{\sqrt{\pi}}{2}\frac{\Gamma\left(\frac{3}{2}\right)\Gamma\left(\frac{3}{2}-\frac{u}{2a}\right)\Gamma\left(\frac{1}{2}-\frac{u}{2a}\right)}{\Gamma(1)\Gamma(1)\Gamma\left(1-\frac{u}{2a}\right)\Gamma\left(1-\frac{u}{2a}\right)}\\
&=\frac{\pi}{8a^2}\frac{\Gamma\left(\frac{1}{2}-\frac{u}{2a}\right)^2}{\Gamma\left(1-\frac{u}{2a}\right)^2}.
\end{align*}
Note that $\zeta(n)=(1-2^{-n})^{-1}t(n)$ and $t(2)=\frac{\pi^2}{8}$, by using \eqref{Eq:Gamma's expansion} and the duplication formula
\begin{align*}
\Gamma\left(\frac{1}{2}-\frac{z}{2}\right)=\frac{\sqrt{\pi}2^z\Gamma(1-z)}{\Gamma\left(1-\frac{z}{2}\right)},
\end{align*}
we have
\begin{align*}
\sum\limits_{k=2}^\infty\left(\sum\limits_{n=1}^{k-1}2^{n-1}G_{0,2a,a}(k,n)\right)u^{k-2}=a^{-2}2^{\frac{2u}{a}}t(2)\exp\left(\sum\limits_{n=2}^\infty\frac{4(2^{n-1}-1)}{na^n(2^n-1)}t(n)u^n\right).
\end{align*}

Similarly, setting $uv=w$ in Theorem \ref{Thm:t-star_N}, we have $\alpha^\star=\frac{a+N-u-v}{N}$, $\beta^\star=\frac{a+N}{N}$. Hence
\begin{align*}
\sum\limits_{k>n\geq1}G_{0,N,a}^\star(k,n)u^{k-n-1}v^{n-1}=\frac{1}{a(a-u-v)}\pFq{3}{2}{\frac{a}{N},\frac{a-u}{N},1}{\frac{a+N}{N},\frac{a+N-u-v}{N}}{1}.
\end{align*}
Setting $v=-2u$ and $N=2a$, then
\begin{align*}
\sum\limits_{k>n\geq1}(-2)^{n-1}G_{0,2a,a}^\star(k,n)u^{k-2}=\frac{1}{a(a+u)}\pFq{3}{2}{\frac{1}{2},\frac{a-u}{2a},1}{\frac{3}{2},\frac{3a+u}{2a}}{1}.
\end{align*}
Therefore we obtain that
\begin{align*}
&\sum\limits_{k=2}^\infty\left(\sum\limits_{n=1}^{k-1}(-2)^{n-1}G_{0,2a,a}^\star(k,n)\right)u^{k-2}\\
&=a^{-2}2^{-\frac{2u}{a}}t(2)\exp\left(\sum\limits_{n=2}^\infty\frac{(-1)^n}{na^n}\frac{4(2^{n-1}-1)}{(2^n-1)}t(n)u^n\right).
\end{align*}
Then the result follows easily.
\qed

Expanding the right-hand side of \eqref{Eq:G_0,2a,a}, we get the following weighted sum formula.
\begin{cor}
For any integer $k\geq2$, we have
\begin{align*}
&\sum\limits_{n=1}^{k-1}2^{n-1}G_{0,2a,a}(k,n)=\sum\limits_{n=1}^{k-1}(-1)^{k-n-1}2^{n-1}G_{0,2a,a}^\star(k,n)\\
&=\sum\limits_{n+n_1+\cdots+n_m=k-2\atop n,m\geq0,n_1,\ldots,n_m\geq2}\frac{2^{n+2m}(2^{n_1-1}-1)\cdots(2^{n_m-1}-1)}{a^kn!m!n_1\cdots n_m(2^{n_1}-1)\cdots(2^{n_m}-1)}t(2)t(n_1)\cdots t(n_m)\log^n(2).
\end{align*}
\end{cor}

\subsection{Proofs of Theorems \ref{Thm:t_N} and \ref{Thm:t-star_N}}\label{Subsec:Proofs}

For an index $\mathbf{k}=(k_1,\ldots,k_n)$, any $a\in\{1,2,\ldots,N\}$ and any nonnegative integer $b$, define
\begin{align}\label{Def:L_N,a,b}
\mathcal{L}_{N,a,b}(\mathbf{k};z)=\sum\limits_{m_1\geq m_2+b\geq\cdots\geq m_n+(n-1)b\geq a+(n-1)b\atop m_i\equiv (n-i)b+a\pmod N}\frac{z^{m_1}}{m_1^{k_1}\cdots m_n^{k_n}}.
\end{align}
Then $\mathcal{L}_{N,a,b}(\mathbf{k};z)$ converges absolutely for $|z|<1$. And we have the following iterated integral representation
\begin{align*}
&\mathcal{L}_{N,a,b}(\mathbf{k};z)=\int_0^z\left(\frac{dt}{t}\right)^{k_1-1}\frac{t^{b-1}dt}{1-t^N}\cdots\left(\frac{dt}{t}\right)^{k_{n-1}-1}\frac{t^{b-1}dt}{1-t^N}\left(\frac{dt}{t}\right)^{k_n-1}\frac{t^{a-1}dt}{1-t^N}.
\end{align*}
For an admissible index $\mathbf{k}$, from the integral representations \eqref{Eq:t_N,a integral} and \eqref{Eq:t_N,a^star integral}, one can easily get that
\begin{align*}
\mathcal{L}_{N,a,N}(\mathbf{k};1)=t_{N,a}(\mathbf{k}),\quad\mathcal{L}_{N,a,0}(\mathbf{k};1)=t_{N,a}^\star(\mathbf{k}).
\end{align*}
By derivation, it is easy to find that
\begin{align}\label{Eq:L_N,a,b derivation}
\frac{d}{dz}\mathcal{L}_{N,a,b}(\mathbf{k};z)=\begin{cases}
\frac{1}{z}\mathcal{L}_{N,a,b}(k_1-1,k_2,\ldots,k_n;z) &  \text{if}\quad k_1>1, \\
\frac{z^{b-1}}{1-z^N}\mathcal{L}_{N,a,b}(k_2,\ldots,k_n;z) & \text{if}\quad n\geq 2,k_1>1, \\
\frac{z^{a-1}}{1-z^N} & \text{if}\quad n=k_1=1.
\end{cases}
\end{align}

For nonnegative integers $k,n,s$, we denote by $I(k,n,s)$ the set of indices of weight $k$, depth $n$ and height $s$, and define the sums
\begin{align*}
G_{N,a,b}(k,n,s;z)=\sum\limits_{\mathbf{k}\in I(k,n,s)}\mathcal{L}_{N,a,b}(\mathbf{k};z),\quad G_{0,N,a,b}(k,n,s;z)=\sum\limits_{\mathbf{k}\in I_0(k,n,s)}\mathcal{L}_{N,a,b}(\mathbf{k};z).
\end{align*}
If the indices set is empty, the sum is treated as zero. And we also set $G_{N,a,b}(0,0,0;z)=1$.
For integers $k,n,s$, using \eqref{Eq:L_N,a,b derivation}, we have
\begin{itemize}
\item[(1)] if $k\geq n+s$ and $n\geq s\geq1$,
\begin{align}\label{Eq:dG_0}
&\frac{d}{dz}G_{0,N,a,b}(k,n,s;z)=\frac{1}{z}\left[G_{N,a,b}(k-1,n,s-1;z)+G_{0,N,a,b}(k-1,n,s;z)\right.\notag\\
&\qquad\qquad\qquad\left.-G_{0,N,a,b}(k-1,n,s-1;z)\right],
\end{align}
\item[(2)] if $k\geq n+s$, $n\geq s\geq0$ and $n\geq2$,
\begin{align}\label{Eq:d(G-G_0)}
\frac{d}{dz}\left[G_{N,a,b}(k,n,s;z)-G_{0,N,a,b}(k,n,s;z)\right]=\frac{z^{b-1}}{1-z^N}G_{N,a,b}(k-1,n-1,s;z).
\end{align}
\end{itemize}

Now we define the generating functions
\begin{align*}
\Phi_{N,a,b}(z)&=\Phi_{N,a,b}(u,v,w;z)=\sum\limits_{k,n,s\geq0}G_{N,a,b}(k,n,s;z)u^{k-n-s}v^{n-s}w^{s},\\
\Phi_{0,N,a,b}(z)&=\Phi_{0,N,a,b}(u,v,w;z)=\sum\limits_{k,n,s\geq0}G_{0,N,a,b}(k,n,s;z)u^{k-n-s}v^{n-s}w^{s-1}\\
&=\sum\limits_{k\geq n+s\atop n\geq s\geq1}G_{0,N,a,b}(k,n,s;z)u^{k-n-s}v^{n-s}w^{s-1}.
\end{align*}
Using \eqref{Eq:L_N,a,b derivation}, \eqref{Eq:dG_0} and \eqref{Eq:d(G-G_0)}, we get that
\begin{align*}
&\frac{d}{dz}\Phi_{0,N,a,b}(z)=\frac{1}{vz}\left(\Phi_{N,a,b}(z)-1-w\Phi_{0,N,a,b}(z)\right)+\frac{u}{z}\Phi_{0,N,a,b}(z),\\
&\frac{d}{dz}\left(\Phi_{N,a,b}(z)-w\Phi_{0,N,a,b}(z)\right)=\frac{vz^{b-1}}{1-z^N}\left(\Phi_{N,a,b}(z)-1\right)+\frac{vz^{a-1}}{1-z^N}.
\end{align*}
Eliminating $\Phi_{N,a,b}(z)$, we obtain the differential equation satisfied by $\Phi_{0,N,a,b}(z)$.

\begin{prop}
$\Phi_{0,N,a,b}=\Phi_{0,N,a,b}(z)$ satisfies the following differential equation
\begin{align}\label{Eq:Phi_0 differential}
z(1-z^N)\Phi_{0,N,a,b}''+[(1-u)(1-z^N)-vz^b]\Phi_{0,N,a,b}'+(uv-w)z^{b-1}\Phi_{0,N,a,b}=z^{a-1}.
\end{align}
\end{prop}

We prove Theorem \ref{Thm:t_N} and Theorem \ref{Thm:t-star_N} by solving the differential equation \eqref{Eq:Phi_0 differential} in the special cases $b=N$ and $b=0$.


\subsubsection{A proof of Theorem \ref{Thm:t_N}}

Let $b=N$ in \eqref{Def:L_N,a,b}. Set
\begin{align*}
G_{0,N,a}(k,n,s;z)=G_{0,N,a,N}(k,n,s;z),\quad\Phi_{0,N,a}(z)=\Phi_{0,N,a,N}(z).
\end{align*}
Using \eqref{Eq:Phi_0 differential}, we get the following differential equation
\begin{align}\label{Eq:t_N differential equation}
z(1-z^N)\Phi_{0,N,a}''+[(1-u)(1-z^N)-vz^N]\Phi_{0,N,a}'+(uv-w)z^{N-1}\Phi_{0,N,a}=z^{a-1}.
\end{align}

We want to find the unique power series solution $\Phi_{0,N,a}(z)=\sum\limits_{n=1}^\infty p_nz^n$. By \eqref{Eq:t_N differential equation}, we see that $p_1,\ldots,p_{a-1},p_{a+1},\ldots,p_N=0$, $p_a=\frac{1}{a(a-u)}$ and
\begin{align*}
p_{n+N}=\frac{n(n-1)+n(1-u+v)-(uv-w)}{(n+N)(n+N-u)}p_n,\quad n\geq1.
\end{align*}
Hence for any $m\geq1$ with $m\not\equiv a\pmod N$, we have $p_m=0$ and for any $n\geq 1$,
\begin{align*}
p_{a+nN}=\frac{(\alpha+n-1)(\beta+n-1)}{(\frac{a}{N}+n)(\frac{a-u}{N}+n)}p_{a+(n-1)N}=\frac{(\alpha)_n(\beta)_n}{(\frac{a}{N}+1)_n(\frac{a-u}{N}+1)_n}\frac{1}{a(a-u)}.
\end{align*}
Therefor we have the following theorem.

\begin{thm}\label{Thm:Phi_0 representation}
For formal variables $u,v,w$, we have
\begin{align*}
\sum\limits_{k\geq n+s,n\geq s\geq1}G_{0,N,a}(k,n,s;z)u^{k-n-s}v^{n-s}w^{s-1}=\frac{z^a}{a(a-u)}\pFq{3}{2}{\alpha,\beta,1}{\frac{a+N}{N},\frac{a+N-u}{N}}{z^N},
\end{align*}
where $\alpha,\beta$ are determined by $\alpha+\beta=\frac{1}{N}(2a-u+v)$ and $\alpha\beta=\frac{1}{N^2}(a(a-u+v)-uv+w)$.
\end{thm}

Setting $z=1$ in Theorem \ref{Thm:Phi_0 representation}, we get Theorem \ref{Thm:t_N}.


\subsubsection{A proof of Theorem \ref{Thm:t-star_N}}

Let $b=0$ in \eqref{Def:L_N,a,b}. Set
\begin{align*}
G_{0,N,a}^\star(k,n,s;z)=G_{0,N,a,0}(k,n,s;z),\quad\Phi_{0,N,a}^\star(z)=\Phi_{0,N,a,0}(z).
\end{align*}
Using \eqref{Eq:Phi_0 differential}, we get the following differential equation
\begin{align}\label{Eq:t_N^star differential equation}
z^2(1-z^N)(\Phi_{0,N,a}^\star)''+z[(1-u)(1-z^N)-v](\Phi_{0,N,a}^\star)'+(uv-w)\Phi_{0,N,a}^\star=z^{a}.
\end{align}

Assume that $\Phi_{0,N,a}^{\star}(z)=\sum\limits_{n=1}^\infty q_nz^n$. By \eqref{Eq:t_N^star differential equation}, we see that $q_1,\ldots,q_{a-1},q_{a+1},\ldots,q_N=0$, $q_a=\frac{1}{a(a-u-v)+uv-w}$ and
\begin{align*}
q_{n+N}=\frac{n(n-u)}{(n+N)(n+N-u-v)+uv-w}q_n,\quad n\geq1.
\end{align*}
Hence for any $m\geq1$ with $m\not\equiv a\pmod N$, we have $q_m=0$, and for any $n\geq 1$
\begin{align*}
q_{a+nN}=\frac{(n+\frac{a}{N}-1)(n+\frac{a-u}{N}-1)}{(\alpha^{\star}+n-1)(\beta^{\star}+n-1)}q_{a+(n-1)N}=\frac{(\frac{a}{N})_n(\frac{a-u}{N})_n}{(\alpha^{\star})_n(\beta^{\star})_n}\frac{1}{a(a-u-v)+uv-w}.
\end{align*}
Therefor we have the following theorem.

\begin{thm}\label{Thm:Phi_0^star representation}
For formal variables $u,v,w$, we have
\begin{align*}
&\sum\limits_{k\geq n+s,n\geq s\geq1}G_{0,N,a}^\star(k,n,s;z)u^{k-n-s}v^{n-s}w^{s-1}\\
&=\frac{z^a}{a(a-u-v)+uv-w}\pFq{3}{2}{\frac{a}{N},\frac{a-u}{N},1}{\alpha^{\star},\beta^{\star}}{z^N},
\end{align*}
where $\alpha^{\star},\beta^{\star}$ are determined by $\alpha^{\star}+\beta^{\star}=\frac{1}{N}(2a+2N-u-v)$ and $\alpha^{\star}\beta^{\star}=\frac{1}{N^2}((a+N)(a+N-u-v)+uv-w)$.
\end{thm}

Setting $z=1$ in Theorem \ref{Thm:Phi_0^star representation}, we get Theorem \ref{Thm:t-star_N}.

\section{Symmetric sum formulas and Hoffman's restricted sum formulas}\label{Sec:Sum formulas}

Recall the construction of stuffle algebra from M. E. Hoffman \cite{Hoffman97}. See also \cite{Li}. Let $N$ be a fixed positive integer and $a\in\{1,2,\ldots,N\}$. Define $\mathfrak{A}_{N,a}=\mathbb{Q} \left \langle \mathrm{x},\mathrm{y}_a \right \rangle$ be the non-commutative polynomial algebra generated by the alphabet $\{\mathrm{x},\mathrm{y}_a\}$ over the rational field $\mathbb{Q}$. Define two subalgebras $\mathfrak{A}_{N,a}^1=\mathbb{Q}+\mathfrak{A}_{N,a}\mathrm{y}_a$ and $\mathfrak{A}_{N,a}^0=\mathbb{Q}+\mathrm{x}\; \mathfrak{A}_{N,a}\mathrm{y}_a$. For a positive integer $k$, set $\mathrm{z}_{k,a}=\mathrm{x}^{k-1}\mathrm{y}_a$. Define two $\mathbb{Q}$-linear maps $t_{N,a}:\mathfrak{A}_{N,a}^0\longrightarrow\mathbb{R}$ and $t_{N,a}^\star:\mathfrak{A}_{N,a}^0\longrightarrow\mathbb{R}$ by $t_{N,a}(1_\mathrm{w})=t_{N,a}^\star(1_\mathrm{w})=1$ and
\begin{align*}
t_{N,a}(\mathrm{z}_{k_1,a}\cdots \mathrm{z}_{k_n,a})=t_{N,a}(k_1,\ldots,k_n),\quad t_{N,a}^\star(\mathrm{z}_{k_1,a}\cdots \mathrm{z}_{k_n,a})=t_{N,a}^\star(k_1,\ldots,k_n),
\end{align*}
where $1_\mathrm{w}$ is the empty word and $k_1,\ldots,k_n$ are positive integers with $k_1>1$. Let $\gamma_{N,a}$ be the algebra automorphism on $\mathfrak{A}_{N,a}$ characterized by $\gamma_{N,a}(\mathrm{x})=\mathrm{x}$ and $\gamma_{N,a}(\mathrm{y}_a)=\mathrm{x}+\mathrm{y}_a$, and define the $\mathbb{Q}$-linear transformation $S_{N,a}:\mathfrak{A}_{N,a}^1\longrightarrow\mathfrak{A}_{N,a}^1$ by
\begin{align*}
S_{N,a}(1_\mathrm{w})=1_\mathrm{w}\quad\text{and}\quad S_{N,a}(\mathrm{w}\mathrm{y}_{a})=\gamma_{N,a}(\mathrm{w})\mathrm{y}_{a}
\end{align*}
for any word $\mathrm{w}\in \mathfrak{A}_{N,a}$. Then by the integral representations \eqref{Eq:t_N,a integral} and \eqref{Eq:t_N,a^star integral}, we have $t^\star_{N,a}=t_{N,a}\circ S_{N,a}$.

Define two commutative products $\ast$ and $\overline{\ast}$ called the stuffle product and the star-stuffle product respectively on $\mathfrak{A}_{N,a}^1$ by $\mathbb{Q}$-bilinearity and the rules:
\begin{align*}
&1_\mathrm{w}\ast \mathrm{w}=\mathrm{w}\ast1_\mathrm{w}=\mathrm{w},\\
&\mathrm{z}_{k,a}\mathrm{w}_1\ast \mathrm{z}_{l,a}\mathrm{w}_2=\mathrm{z}_{k,a}(\mathrm{w}_1\ast \mathrm{z}_{l,a}\mathrm{w}_2)+\mathrm{z}_{l,a}(\mathrm{z}_{k,a}\mathrm{w}_1\ast \mathrm{w}_2)+\mathrm{z}_{k+l,a}(\mathrm{w}_1\ast \mathrm{w}_2),\\
&1_\mathrm{w}\overline{\ast} \mathrm{w}=\mathrm{w}\overline{\ast}1_\mathrm{w}=\mathrm{w},\\
&\mathrm{z}_{k,a}\mathrm{w}_1\overline{\ast} \mathrm{z}_{l,a}\mathrm{w}_2=\mathrm{z}_{k,a}(\mathrm{w}_1\overline{\ast} \mathrm{z}_{l,a}\mathrm{w}_2)+\mathrm{z}_{l,a}(\mathrm{z}_{k,a}\mathrm{w}_1\overline{\ast}\mathrm{w}_2)-\mathrm{z}_{k+l,a}(\mathrm{w}_1\overline{\ast}\mathrm{w}_2),
\end{align*}
where $k,l$ are positive integers and $\mathrm{w},\mathrm{w}_1,\mathrm{w}_2$ are words in $\mathfrak{A}_{N,a}^1$. Then we get commutative algebras $\mathfrak{A}_{N,a,\ast}^1$ and $\mathfrak{A}_{N,a,\overline{\ast}}^1$, their subalgebras $\mathfrak{A}_{N,a,\ast}^0$ and $\mathfrak{A}_{N,a,\overline{\ast}}^0$. It is easy to see that the maps $t_{N,a}:\mathfrak{A}_{N,a,\ast}^0\longrightarrow\mathbb{R}$ and $t_{N,a}^\star:\mathfrak{A}_{N,a,\overline{\ast}}^0\longrightarrow\mathbb{R}$ are algebra homomorphisms.


Denote by $\mathcal{\mathcal{P}}_n$ the set of all partitions of the set $\{1,2,\dots,n\}$. For a partition $\Pi=\{P_1,P_2,\ldots,P_l\}\in \mathcal{\mathcal{P}}_n$, let
\begin{align*}
c(\Pi)=\prod\limits_{j=1}^l(|P_j|-1)!,\quad\widetilde{c}(\Pi)=(-1)^{n-l}c(\Pi),
\end{align*}
and for any maximal height index $\mathbf{k}=(k_1,\ldots,k_n)$, define
\begin{align*}
t_{N,a}(\mathbf{k},\Pi)=\prod\limits_{j=1}^lt_{N,a}\left(\sum\limits_{l\in P_j}k_l\right),\quad t_{N,a}^\star(\mathbf{k},\Pi)=\prod\limits_{j=1}^lt_{N,a}^\star\left(\sum\limits_{l\in P_j}k_l\right).
\end{align*}
Then applying the maps $t_{N,a}$ and $t^{\star}_{N,a}$ to \cite[Lemma 5.1]{Li-Qin}, we get the symmetric sum formulas for multiple $t$-(star) values of level $N$, which generalizes the results for multiple zeta(-star) values \cite[Theorems 2.1 and 2.2]{Hoffman92} and multiple $t$-(star) values \cite[Theorems 3.2 and 3.5]{Hoffman19}.

\begin{prop}\label{Prop:SymSum}
For an index $\mathbf{k}=(k_1,\ldots,k_n)$ with $k_1,\ldots,k_n>1$, we have
\begin{align*}
\sum\limits_{\sigma\in S_n}t_{N,a}(k_{\sigma(1)},\ldots,k_{\sigma(n)})=\sum\limits_{\Pi\in\mathcal{\mathcal{P}}_n}\widetilde{c}(\Pi)t_{N,a}(\mathbf{k},\Pi),\\
\sum\limits_{\sigma\in S_n}t_{N,a}^\star(k_{\sigma(1)},\ldots,k_{\sigma(n)})=\sum\limits_{\Pi\in\mathcal{\mathcal{P}}_n}c(\Pi)t_{N,a}^\star(\mathbf{k},\Pi),
\end{align*}
where $S_n$ is the symmetric group of degree $n$.
\end{prop}

Taking $k_1=\cdots=k_n$ in this result, we get a generating function for multiple $t$-(star) values of level $N$ with repeated arguments.

\begin{cor}\label{Cor:Rep-Arg}
For any integer $k>1$, we have
\begin{align}
&1+\sum\limits_{n=1}^\infty t_{N,a}(\{k\}^n)x^{kn}=\exp\left(\sum\limits_{n=1}^\infty \frac{(-1)^{n-1}t_{N,a}(kn)x^{kn}}{n}\right),\label{Eq:t-N,a-generating-function}\\
&1+\sum\limits_{n=1}^\infty t_{N,a}^\star(\{k\}^n)x^{kn}=\exp\left(\sum\limits_{n=1}^\infty \frac{t_{N,a}(kn)x^{kn}}{n}\right).\label{Eq:t-star-N,a-generating-function}
\end{align}
\end{cor}

\proof
We prove \eqref{Eq:t-star-N,a-generating-function} and one can prove \eqref{Eq:t-N,a-generating-function} similarly. Using Proposition \ref{Prop:SymSum} with $k_1=\cdots=k_n=k$, we get
\begin{align*}
t_{N,a}^\star(\{k\}^n)&=\frac{1}{n!}\sum\limits_{l=1}^n\sum\limits_{\{P_1,\ldots,P_l\}\in\mathcal{P}_n}\prod\limits_{s=1}^l(p_s-1)!t_{N,a}^\star(kp_s),
\end{align*}
where $p_s=|P_s|$. For fixed $p_1,\ldots,p_l$, we have
\begin{align*}
\left|\{\{P_1,\ldots,P_l\}\in \mathcal{P}_n\mid p_s=|P_s|, s=1,\ldots,l\}\right|=\frac{1}{m_1!m_2!\cdots m_n!}\frac{n!}{p_1!p_2!\cdots p_l!},
\end{align*}
where $m_i=\left|\{s\mid p_s=i\}\right|$. Then we obtain
\begin{align*}
&t_{N,a}^\star(\{k\}^n)\\
&=\sum\limits_{l=1}^n\sum\limits_{p_1+\cdots+p_l=n\atop p_1,\ldots,p_l>0}\frac{1}{m_1!\cdots m_n!p_1\cdots p_l}t_{N,a}^\star(kp_1)\cdots t_{N,a}^\star(kp_l)\\
&=\sum\limits_{m_1+2m_2+\cdots+nm_n=n\atop m_1,\ldots,m_n\geq 0}\frac{1}{m_1!m_2!\cdots m_n!}\left(t_{N,a}^\star(k)\right)^{m_1}\left(\frac{t_{N,a}^\star(2k)}{2}\right)^{m_2}\cdots\left(\frac{t_{N,a}^\star(nk)}{n}\right)^{m_n}.
\end{align*}
Hence we easily get the desired generating function.
\qed

Setting $N=a=1$ in Corollary \ref{Cor:Rep-Arg}, we get the generating function for multiple zeta-(star) values with repeated arguments \cite[Proposition 3]{Ihara-Kajikawa-Ohno-Okuda}. And setting $N=2$ and $a=1$ in Corollary \ref{Cor:Rep-Arg}, we obtain the generating function for multiple $t$-(star) values with repeated arguments. Note that using $t(k)=(1-2^{-k})\zeta(k)$, M. E. Hoffman \cite[Theorems 3.4 and 3.6]{Hoffman19} in fact showed that
\begin{align}
&1+\sum\limits_{n=1}^\infty t(\{k\}^n)x^{kn}=\frac{Z_k(x)}{Z_k\left(\frac{x}{2}\right)},\notag\\
&1+\sum\limits_{n=1}^\infty t^\star(\{k\}^n)x^{kn}=\frac{Z_k(e^{\frac{\pi i}{k}}\frac{x}{2})}{Z_k(e^{\frac{\pi i}{k}}x)},\label{Eq:t-star-generating-function}
\end{align}
where $Z_k(x)=1+\sum\limits_{n=1}^\infty \zeta(\{k\}^n)x^{kn}$. We need the following generating function for $t^\star(\{2\}^n)$ which is important for the evaluations in Section \ref{Sec:Evaluations}.

\begin{prop}
We have
\begin{align}\label{Eq:generating function of t-star}
1+\sum\limits_{n=1}^\infty t^\star(\{2\}^n)x^{2n}=\frac{1}{\cos(\frac{\pi x}{2})}.
\end{align}
\end{prop}

\proof
Using \eqref{Eq:t-star-generating-function} and \cite[Eq. (34)]{Broadhurst-Borwein-Bradley}
\begin{align*}
Z_{2l}(x)=\frac{1}{(i\pi x)^l}\prod\limits_{j=1}^l\sin\left(e^{\frac{(2j-1)\pi i}{2l}}\pi x\right),
\end{align*}
we easily get the result.
\qed

\begin{rem}
Since $\sum\limits_{a=1}^Nt_{N,a}(k)=\zeta(k)$, for any integer $k>1$, we have
\begin{align*}
\prod\limits_{a=1}^N\left(1+\sum\limits_{n=1}^\infty t_{N,a}(\{k\}^n)x^{kn}\right)=Z_k(x).
\end{align*}
\end{rem}

Applying the maps $t_{N,a}$ and $t^{\star}_{N,a}$ to \cite[Proposition 3.27]{Li-Qin2021}, we get Hoffman's restricted sum formula for multiple $t$-(star) values of level $N$, which generalizes the results for multiple zeta values \cite[Theorem A]{Chen-Chung-Eie}, multiple zeta-star values \cite[Theorem 1.1]{Chen-Chung-Eie2} and multiple $t$-values \cite[Theorem 1]{Shen-Jia}.
\begin{prop}
For any integer $m>1$, we have
\begin{align*}
&\sum\limits_{k_1+\cdots+k_n=k\atop k_i\geq1}t_{N,a}(mk_1,\ldots,mk_n)=\sum\limits_{j=0}^{k-n}(-1)^{k-n-j}\binom{k-j}{n}t_{N,a}^\star(\{m\}^j)t_{N,a}(\{m\}^{k-j}),\\
&\sum\limits_{k_1+\cdots+k_n=k\atop k_i\geq1}t_{N,a}^\star(mk_1,\ldots,mk_n)=\sum\limits_{j=0}^{k-n}(-1)^{j}\binom{k-j}{n}t_{N,a}(\{m\}^j)t_{N,a}^\star(\{m\}^{k-j}).
\end{align*}
\end{prop}

\section{Evaluations with one-two-three indices}\label{Sec:Evaluations}

In this section, we consider multiple $t$-(star) values of level $2$. For simplicity, we denote $\mathfrak{A}_{2,1}^1$ and $S_{2,1}$ defined in Section \ref{Sec:Sum formulas} by $\mathfrak{A}^1$ and $S$, respectively. And for a positive integer $k$, set $\mathrm{z}_k=\mathrm{z}_{k,1}\in\mathfrak{A}^1$. Recall that for an index $\mathbf{k}$, there is a stuffle regularization $t_{\ast}^{V}(\mathbf{k})$ of multiple $t$-values with regularization parameter $t_{\ast}^{V}(1)=V$ in \cite{Charlton}. For a word $\mathrm{w}=\mathrm{z}_{k_1}\cdots \mathrm{z}_{k_n}\in\mathfrak{A}^1$, define $t_{\ast}^{V}(\mathrm{w})=t_{\ast}^{V}(k_1,\ldots,k_n)$. Then the stuffle regularization map $t_{\ast}^V$ can be extended to $\mathfrak{A}^1$ by $\mathbb{Q}$-linearities. The star-stuffle regularization map $t^{\star,V}_{\overline{\ast}}$ on $\mathfrak{A}^1$ is defined by
\begin{align*}
t^{\star,V}_{\overline{\ast}}(\mathrm{w})
=t_{\ast}^{V}(S(\mathrm{w})),\quad \text{for any}\quad \mathrm{w}\in\mathfrak{A}^1.
\end{align*}
Similarly, for an index $\mathbf{k}=(k_1,\ldots,k_n)$, set $t^{\star,V}_{\overline{\ast}}(\mathbf{k})=t^{\star,V}_{\overline{\ast}}(\mathrm{z}_{k_1}\cdots \mathrm{z}_{k_n})$. Then we have $t^{\star,V}_{\overline{\ast}}(1)=t^{\star,V}_{\overline{\ast}}(\mathrm{y}_1)=V$. 

In this section, we get two evaluations for the multiple $t$-star values $t^\star(\{2\}^a,3,\{2\}^b)$ and $t^{\star,V}_{\overline{\ast}}(\{2\}^a,1,\{2\}^b)$, which are analogous to the evaluations of the multiple $t$-values established in \cite{Murakami2021} and \cite{Charlton} respectively.

\begin{thm}\label{Thm:t_star two-three}
For any integers $a,b\geq 0$, we have
\begin{align}\label{Eq:t_star two-three}
&t^\star(\{2\}^a,3,\{2\}^b)\notag\\
=&\sum\limits_{r=1}^{a+b+1}2^{-2r}\left[(1-2^{-2r})\binom{2r}{2a+1}+\binom{2r}{2b+1}\right]t^\star\left(\{2\}^{a+b+1-r}\right)\zeta(2r+1).
\end{align}
\end{thm}

\begin{thm}\label{Thm:t_star one-two}
For any integers $a,b\geq 0$, we have
\begin{align}\label{Eq:t_star one-two}
t^{\star,V}_{\overline{\ast}}(\{2\}^a,1,\{2\}^b)&=\sum\limits_{r=1}^{a+b}2^{-2r}\left[\binom{2r}{2a}+(1-2^{-2r})\binom{2r}{2b}\right]t^\star\left(\{2\}^{a+b-r}\right)\zeta(2r+1)\notag\\
&\qquad\qquad\qquad+\delta_{a=0}(V-\log (2))t^\star(\{2\}^b)+\delta_{b=0}\log(2)t^\star(\{2\}^a).
\end{align}
where $\delta$ is the Kronecker symbol. 
\end{thm}

We will prove Theorem \ref{Thm:t_star two-three} and Theorem \ref{Thm:t_star one-two} by using the generating functions. Furthermore, we will show the algebraic equivalences between the above evaluation formulas and the corresponding evaluations for multiple $t$-values.

\subsection{Evaluation of $t^\star(\{2\}^a,3,\{2\}^b)$}

To evaluate $t^\star(\{2\}^a,3,\{2\}^b)$, we need Murakami's generating function for $t(\{2\}^a,3,\{2\}^b)$ in \cite{Murakami2021}. For an admissible index $(k_1,\ldots,k_n)$, set $\widetilde{t}(k_1,\ldots,k_n)=2^{k_1+\cdots+k_n}t(k_1,\ldots,k_n)$. Let $K(a,b)=\widetilde{t}(\{2\}^a,3,\{2\}^b)$ and $G(x,y)$ be the generating function
$$G(x,y)=\sum\limits_{a,b\geq0}(-1)^{a+b}K(a,b)x^{2a+1}y^{2b+1}.$$
By the definition of $K(a,b)$, $G(x,y)$ can be written as
\begin{align}
G(x,y)=xy\sum\limits_{m=1}^\infty\prod\limits_{k=m+1}^\infty\left(1-\frac{x^2}{\left(k-\frac{1}{2}\right)^2}\right)\frac{1}{\left(m-\frac{1}{2}\right)^3}\prod\limits_{l=1}^{m-1}\left(1-\frac{y^2}{\left(l-\frac{1}{2}\right)^2}\right).\label{Eq:G(x,y) def}
\end{align}
From \cite[Proposition 13]{Murakami2021}, $G(x,y)$ is displayed as
\begin{align}\label{Eq:G(x,y)}
G(x,y)=\cos\pi x[A(x+y)-A(x-y)]+\cos\pi y[B(x+y)-B(x-y)],
\end{align}
where $A(z)$ and $B(z)$ are the power series
\begin{align*}
A(z)=\sum\limits_{r=1}^\infty\zeta(2r+1)z^{2r},\quad B(z)=\sum\limits_{r=1}^\infty(1-2^{-2r})\zeta(2r+1)z^{2r}.
\end{align*}

\begin{proofof}{Theorem \ref{Thm:t_star two-three}}
Let $K^\star(a,b)=\widetilde{t}^\star(\{2\}^a,3,\{2\}^b)$ and $G^\star(x,y)$ be the generating function
\begin{align*}
G^\star(x,y)&=\sum\limits_{a,b\geq0}K^\star(a,b)x^{2a}y^{2b},
\end{align*}
where $\widetilde{t}^{\star}(k_1,\ldots,k_n)=2^{k_1+\cdots+k_n}t^{\star}(k_1,\ldots,k_n)$. By the definitions of $K^\star(a,b)$, we get
\begin{align*}
G^\star(x,y)&=\sum\limits_{m=1}^\infty\prod\limits_{k=m}^\infty\left(1-\frac{x^2}{\left(k-\frac{1}{2}\right)^2}\right)^{-1}\frac{1}{\left(m-\frac{1}{2}\right)^3}\prod\limits_{l=1}^m\left(1-\frac{y^2}{\left(l-\frac{1}{2}\right)^2}\right)^{-1}.\\
\end{align*}
Using \eqref{Eq:G(x,y) def} together with \eqref{Eq:G(x,y)}, we have
\begin{align*}
G^\star(x,y)&=\frac{G(y,x)}{xy\cos\pi x\cos\pi y}\\
&=\frac{A(x+y)-A(x-y)}{xy\cos\pi x}+\frac{B(x+y)-B(x-y)}{xy\cos\pi y}.
\end{align*}
By \eqref{Eq:generating function of t-star} and the identity
\begin{align}
2\sum\limits_{n=0}^{r-1}\binom{2r}{2n+1}x^{2r-2n-1}y^{2n+1}=(x+y)^{2r}-(x-y)^{2r},\label{Eq:binom{2r}{2n+1}}
\end{align}
we find that
\begin{align*}
G^\star(x,y)&=2\sum\limits_{n=0}^\infty\sum\limits_{r=1}^\infty\sum\limits_{k=0}^{r-1}\binom{2r}{2k+1}\widetilde{t}^\star(\{2\}^n)\zeta(2r+1)x^{2n+2r-2k-2}y^{2k}\\
&\qquad\qquad+2\sum\limits_{n=0}^\infty\sum\limits_{r=1}^\infty\sum\limits_{l=0}^{r-1}(1-2^{-2r})\binom{2r}{2l+1}\widetilde{t}^\star(\{2\}^n)\zeta(2r+1)x^{2r-2l-2}y^{2n+2l}.
\end{align*}
Hence we get the desired result by comparing the coefficient of $x^{2a}y^{2b}$ on both sides.
\end{proofof}

The evaluation of $t(\{2\}^a,3,\{2\}^b)$ in \cite[Theorem 3]{Murakami2021} shows that
\begin{align}\label{Eq:t-two-three}
&t(\{2\}^a,3,\{2\}^b)\notag\\
=&\sum\limits_{r=1}^{a+b+1}(-1)^{r-1}2^{-2r}\left[(1-2^{-2r})\binom{2r}{2a+1}
+\binom{2r}{2b+1}\right]t\left(\{2\}^{a+b+1-r}\right)\zeta(2r+1).
\end{align}
We prove the  equivalence of \eqref{Eq:t_star two-three} and \eqref{Eq:t-two-three}.

\begin{prop}
\eqref{Eq:t_star two-three} holds for any integers $a,b\geq 0$ if and only if \eqref{Eq:t-two-three} holds for any integers $a,b\geq 0$.
\end{prop}

\proof
We give an algebraic proof here. Using \cite[Eq. (2.4)]{Li}, we get that
\begin{align}\label{Eq:T-star vs T}
\mathcal{T}^\star(u,v)=\mathcal{T}(v,u)\ast \mathcal{S}^\star(u)\ast \mathcal{S}^\star(v),
\end{align}
where
\begin{align*}
\mathcal{T}(u,v)=\sum\limits_{a,b=0}^\infty(-1)^{a+b}\mathrm{z}_2^a\mathrm{z}_3\mathrm{z}_2^bu^{2a}v^{2b},\quad \mathcal{T}^\star(u,v)=\sum\limits_{a,b=0}^\infty S(\mathrm{z}_2^a\mathrm{z}_3\mathrm{z}_2^b)u^{2a}v^{2b}
\end{align*}
and
\begin{align*}
\mathcal{S}^\star(u)=\sum\limits_{n=0}^\infty S(\mathrm{z}_2^n)u^{2n}.
\end{align*}

On the other hand, let
\begin{align*}
\hat{\mathcal{T}}(u,v)=\sum\limits_{a,b=0}^\infty(-1)^{a+b}\hat{\mathcal{H}}(a,b)u^{2a}v^{2b},\quad\hat{\mathcal{T}}^\star(u,v)
=\sum\limits_{a,b=0}^\infty\hat{\mathcal{H}}^\star(a,b)u^{2a}v^{2b},
\end{align*}
where
\begin{align*}
&\hat{\mathcal{H}}(a,b)=\sum\limits_{r=1}^{a+b+1}(-1)^{r-1}2^{-2r}\left[(1-2^{-2r})\binom{2r}{2a+1}+\binom{2r}{2b+1}\right]\mathrm{z}_{2r+1}\ast \mathrm{z}_2^{a+b+1-r},\\
&\hat{\mathcal{H}}^\star(a,b)=\sum\limits_{r=1}^{a+b+1}2^{-2r}\left[(1-2^{-2r})\binom{2r}{2a+1}+\binom{2r}{2b+1}\right]\mathrm{z}_{2r+1}\ast S(\mathrm{z}_2^{a+b+1-r}).
\end{align*}
Set
\begin{align*}
\mathcal{A}(u)=\sum\limits_{r=1}^\infty \mathrm{z}_{2r+1}u^{2r},\quad
\mathcal{B}(u)=\sum\limits_{r=1}^\infty(1-2^{-2r})\mathrm{z}_{2r+1}u^{2r}
\end{align*}
and
\begin{align*}
\mathcal{S}(u)=\sum\limits_{n=0}^\infty (-1)^n\mathrm{z}_2^nu^{2n}.
\end{align*}
Then we find that
\begin{align}
\hat{\mathcal{T}}(u,v)=&\frac{\mathcal{S}(u)}{2uv}\ast\left[\mathcal{A}\left(\frac{u+v}{2}\right)-\mathcal{A}\left(\frac{u-v}{2}\right)\right]\notag\\
&+\frac{\mathcal{S}(v)}{2uv}\ast\left[\mathcal{B}\left(\frac{u+v}{2}\right)-\mathcal{B}\left(\frac{u-v}{2}\right)\right],\label{Eq:hat-T}\\
\hat{\mathcal{T}}^\star(v,u)=&\frac{\mathcal{S}^\star(u)}{2uv}\ast\left[\mathcal{A}\left(\frac{u+v}{2}\right)-\mathcal{A}\left(\frac{u-v}{2}\right)\right]\notag\\
&+\frac{\mathcal{S}^\star(v)}{2uv}\ast\left[\mathcal{B}\left(\frac{u+v}{2}\right)-\mathcal{B}\left(\frac{u-v}{2}\right)\right].\label{Eq:hat-T-star}
\end{align}
We give the proof of \eqref{Eq:hat-T} and one can prove \eqref{Eq:hat-T-star} similarly. In fact, we have
\begin{align*}
\hat{\mathcal{T}}(u,v)&=\sum\limits_{a,b=0}^\infty\sum\limits_{r=1}^\infty(-1)^{a+b+r-1}2^{-2r}\left[(1-2^{-2r})\binom{2r}{2a+1}+\binom{2r}{2b+1}\right]\\
&\qquad\times \mathrm{z}_{2r+1}\ast \mathrm{z}_2^{a+b+1-r}u^{2a}v^{2b}\\
&=\sum\limits_{r\geq 1,s\geq0}\sum\limits_{a=0}^{r-1}(-1)^s2^{-2r}(1-2^{-2r})\binom{2r}{2a+1}\mathrm{z}_{2r+1}\ast \mathrm{z}_2^{s}u^{2a}v^{2s+2r-2b-2}\\
&\qquad+\sum\limits_{r\geq 1,s\geq0}\sum\limits_{b=0}^{r-1}(-1)^s2^{-2r}\binom{2r}{2b+1}\mathrm{z}_{2r+1}\ast \mathrm{z}_2^{s}u^{2s+2r-2a-2}v^{2b}\\
&=\sum\limits_{s=0}^\infty(-1)^s\mathrm{z}_2^sv^{2s}\ast\sum\limits_{r=1}^\infty2^{-2r}(1-2^{-2r})\mathrm{z}_{2r+1}\sum\limits_{a=0}^{r-1}\binom{2r}{2a+1}u^{2a}v^{2r-2b-2}\\
&\qquad+\sum\limits_{s=0}^\infty(-1)^s\mathrm{z}_2^su^{2s}\ast\sum\limits_{r=1}^\infty2^{-2r}\mathrm{z}_{2r+1}\sum\limits_{b=0}^{r-1}\binom{2r}{2b+1}u^{2r-2a-2}v^{2b}.
\end{align*}
Hence by \eqref{Eq:binom{2r}{2n+1}} and the definitions of $\mathcal{A}(u),\mathcal{B}(u)$ and $\mathcal{S}(u)$, we get \eqref{Eq:hat-T} easily.

Notice the fact $\mathcal{S}(u)\ast \mathcal{S}^\star(u)=1$ (see \cite[Corollary 1]{Ihara-Kajikawa-Ohno-Okuda}). Then we have
\begin{align}\label{Eq:hat-T-star vs hat-T}
\hat{\mathcal{T}}^\star(u,v)=\hat{\mathcal{T}}(v,u)\ast \mathcal{S}^\star(u)\ast \mathcal{S}^\star(v).
\end{align}
Hence, \eqref{Eq:T-star vs T} and \eqref{Eq:hat-T-star vs hat-T} imply the equivalence of the evaluation formulas \eqref{Eq:t_star two-three} and \eqref{Eq:t-two-three}.\qed

\subsection{Evaluation of $t^{\star,V}_{\overline{\ast}}(\{1\}^a,2,\{1\}^b)$}

Analogous to \cite[Definition 3.4]{Charlton}, for $s_1,\ldots,s_n\in\mathbb{Z}_{\geq 0}$, the multiple $t$-star polylogarithm is defined by
\begin{align*}
\Ti^\star_{s_1,\ldots,s_n}(z_1,\ldots,z_n)=\sum\limits_{m_1\geq\cdots\geq m_n>0}\frac{z_1^{2m_1-1}\cdots z_n^{2m_n-1}}{(2m_1-1)^{s_1}\cdots(2m_n-1)^{s_n}},
\end{align*}
which converges when $|z_1\cdots z_i|<1$ for $i=1,\ldots,n$. Note that for $s_1>1$ and $z_1=\cdots=z_n=1$, $\Ti^\star_{s_1,\ldots,s_n}(1,\ldots,1)$ is exactly the multiple $t$-star value $t^\star(s_1,\ldots,s_n)$. When $n=s_1=1$, we get
\begin{align}\label{Eq:Ti-star-1-z}
\Ti^\star_1(z)=\sum\limits_{m=1}^\infty\frac{z^{2m-1}}{2m-1}=\tanh^{-1}(z),
\end{align}
which tends to infinity when $z\to 1^{-}$.

We calculate the generating function
\begin{align}\label{Eq:generation functio of Ti-star-2-a-1-2-b}
&\sum\limits_{a,b\geq0}\Ti^\star_{\{2\}^a,1,\{2\}^b}(\{1\}^a,z,\{1\}^b)(2x)^{2a}(2y)^{2b}\notag\\
&=\sum\limits_{r=1}^\infty\prod\limits_{k=r}^\infty\left(1-\frac{(2x)^2}{(2k-1)^2}\right)^{-1}\frac{z^{2r-1}}{2r-1}
\prod\limits_{l=1}^r\left(1-\frac{(2y)^2}{(2l-1)^2}\right)^{-1}\notag\\
&=\frac{1}{\cos\pi x}\sum\limits_{r=1}^\infty\prod_{k=1}^{r-1}\frac{(2k-1)^2-(2x)^2}{(2k-1)^2}\frac{z^{2r-1}}{2r-1}
\prod\limits_{l=1}^r\frac{(2l-1)^2}{(2l-1)^2-(2y)^2}\notag\\
&=\frac{1}{(1-4y^2)\cos\pi x}\sum\limits_{r=1}^\infty\frac{\left(\frac{1}{2}+x\right)\cdots\left(r-\frac{3}{2}+x\right)\left(\frac{1}{2}-x\right)\cdots\left(r-\frac{3}{2}-x\right)(2r-1)z^{2r-1}}{\left(\frac{3}{2}+y\right)\cdots\left(r-\frac{1}{2}+y\right)\left(\frac{3}{2}-y\right)\cdots\left(r-\frac{1}{2}-y\right)}\notag\\
&=\frac{z}{(1-4y^2)\cos\pi x}\pFq{4}{3}{1,\frac{3}{2},\frac{1}{2}-x,\frac{1}{2}+x}{\frac{1}{2},\frac{3}{2}-y,\frac{3}{2}+y}{z^2}.
\end{align}
By the stuffle product, we have
\begin{align*}
\Ti^\star_{1,\{2\}^b}(z,\{1\}^b)&=\Ti^\star_{1}(z)\Ti^\star_{\{2\}^{b}}(\{1\}^{b})-\sum\limits_{i=1}^b\Ti^\star_{\{2\}^i,1,\{2\}^{b-i}}(\{1\}^i,z,\{1\}^{b-i})\notag\\
&\qquad+\sum\limits_{j=0}^{b-1}\Ti^\star_{\{2\}^j,3,\{2\}^{b-1-j}}(\{1\}^j,z,\{1\}^{b-1-j}).
\end{align*}
Hence
\begin{align*}
t^{\star,V=0}_{\overline{\ast}}(1,\{2\}^b)=-\sum\limits_{i=1}^bt^\star(\{2\}^i,1,\{2\}^{b-i})+\sum\limits_{j=0}^{b-1}t^\star(\{2\}^j,3,\{2\}^{b-1-j}).
\end{align*}
Then by \eqref{Eq:generating function of t-star} and \eqref{Eq:Ti-star-1-z}, the divergent part (as $z\to 1^-$) of the generating function can be written as
\begin{align}\label{Eq:generation functio of Ti-star-1-2-b}
\sum\limits_{b=0}^\infty\Ti^\star_{1,\{2\}^b}(z,\{1\}^b)(2y)^{2b}
&=\sum\limits_{b=0}^\infty\Ti^\star_{1}(z)t^\star(\{2\}^b)(2y)^{2b}
+\sum\limits_{b=0}^\infty t^{\star,V=0}_{\overline{\ast}}(1,\{2\}^b)(2y)^{2b}\notag\\
&=\frac{\tanh^{-1}(z)}{\cos\pi y}+\sum\limits_{b=0}^\infty t^{\star,V=0}_{\overline{\ast}}(1,\{2\}^b)(2y)^{2b}.
\end{align}
Subtracting \eqref{Eq:generation functio of Ti-star-1-2-b} from \eqref{Eq:generation functio of Ti-star-2-a-1-2-b} and letting $z\to1^{-}$, we get the generating function of the star-stuffle regularized (at $V=0$) multiple $t$-star values 
\begin{align}
&\sum\limits_{a,b\geq 0}t^{\star,V=0}_{\overline{\ast}}(\{2\}^a,1,\{2\}^b)(2x)^{2a}(2y)^{2b}\notag\\
=&\lim_{z\to 1^-}\left\{\frac{z}{(1-4y^2)\cos\pi x}\pFq{4}{3}{1,\frac{3}{2},\frac{1}{2}-x,\frac{1}{2}+x}{\frac{1}{2},\frac{3}{2}-y,\frac{3}{2}+y}{z^2}-\frac{\tanh^{-1}(z)}{\cos\pi y}\right\}\notag\\
=&\frac{1}{2\cos\pi y}[A(x+y)+A(x-y)-2\log(2)]\notag\\
&\qquad+\frac{1}{2\cos\pi x}[B(x+y)+B(x-y)+2\log(2)],\label{Eq:generation function at V^star=0}
\end{align}
where the last step follows from \cite[Eqs. (14), (17)]{Charlton}.

Then by comparing the coefficient of $x^{2a}y^{2b}$ on both sides of \eqref{Eq:generation function at V^star=0}, we have
\begin{align}\label{Eq:generation function at V^star}
t^{\star,V=0}_{\overline{\ast}}(\{2\}^a,1,\{2\}^b)&=\sum\limits_{r=1}^{a+b}2^{-2r}\left[\binom{2r}{2a}+(1-2^{-2r})\binom{2r}{2b}\right]t^\star\left(\{2\}^{a+b-r}\right)\zeta(2r+1)\notag\\
&\qquad\qquad-\delta_{a=0}\log(2) t^\star(\{2\}^b)+\delta_{b=0}\log(2)t^\star(\{2\}^a).
\end{align}
Notice that
\begin{align*}
t^{\star,V}_{\overline{\ast}}(1,\{2\}^b)=V t^\star(\{2\}^b)+t^{\star,V=0}_{\overline{\ast}}(1,\{2\}^b),
\end{align*}
we complete the proof of Theorem \ref{Thm:t_star one-two} by giving the necessary correction terms to add to the right-hand side of \eqref{Eq:generation function at V^star}.
\qed

The evaluation of $t_{\ast}^V(\{2\}^a,1,\{2\}^b)$ in \cite[Theorem 1.1]{Charlton} shows that
\begin{align}\label{Eq:t-one-two}
t^V_{\ast}(\{2\}^a,1,\{2\}^b)&=\sum\limits_{r=1}^{a+b}(-1)^r2^{-2r}\left[\binom{2r}{2a}+(1-2^{-2r})\binom{2r}{2b}\right]t\left(\{2\}^{a+b-r}\right)\zeta(2r+1)\notag\\
&\qquad\qquad\qquad+\delta_{a=0}(V-\log(2))t(\{2\}^b)+\delta_{b=0}\log(2)t(\{2\}^a).
\end{align}
We give an algebraic proof of the equivalence of \eqref{Eq:t_star one-two} and \eqref{Eq:t-one-two}.

\begin{prop}
\eqref{Eq:t_star one-two} holds for any integers $a,b\geq 0$ if and only if  \eqref{Eq:t-one-two} holds for any integers $a,b\geq 0$.
\end{prop}

\proof
Using \cite[Eq. (2.4)]{Li}, we get that
\begin{align}\label{Eq:P vs P-star}
\mathcal{P}^\star(u,v)=\mathcal{P}(v,u)\ast \mathcal{S}^\star(u)\ast \mathcal{S}^\star(v),
\end{align}
where
\begin{align*}
\mathcal{P}(u,v)=\sum\limits_{a,b=0}^\infty(-1)^{a+b}\mathrm{z}_2^a\mathrm{z}_1\mathrm{z}_2^bu^{2a}v^{2b},\quad \mathcal{P}^\star(u,v)=\sum\limits_{a,b=0}^\infty S(\mathrm{z}_2^a\mathrm{z}_1\mathrm{z}_2^b)u^{2a}v^{2b}.
\end{align*}

Let
\begin{align*}
\hat{\mathcal{P}}(u,v)=\sum\limits_{a,b=0}^\infty(-1)^{a+b}\hat{\mathcal{J}}(a,b)u^{2a}v^{2b},\quad\hat{\mathcal{P}}^\star(u,v)
=\sum\limits_{a,b=0}^\infty\hat{\mathcal{J}}^\star(a,b)u^{2a}v^{2b},
\end{align*}
where
\begin{align*}
&\hat{\mathcal{J}}(a,b)=\sum\limits_{r=1}^{a+b}(-1)^r2^{-2r}\left[\binom{2r}{2a}+(1-2^{-2r})\binom{2r}{2b}\right]\mathrm{z}_{2r+1}\ast \mathrm{z}_2^{a+b-r}\\
&\qquad\qquad\qquad\qquad\qquad\qquad\qquad+\delta_{a=0}(\mathrm{z}_1-\log(2))\ast \mathrm{z}_2^b+\delta_{b=0}\log(2)\mathrm{z}_2^a,\\
&\hat{\mathcal{J}}^\star(a,b)=\sum\limits_{r=1}^{a+b}2^{-2r}\left[\binom{2r}{2a}+(1-2^{-2r})\binom{2r}{2b}\right]\mathrm{z}_{2r+1}\ast S(\mathrm{z}_2^{a+b-r})\\
&\qquad\qquad\qquad\qquad\qquad\qquad+\delta_{a=0}(\mathrm{z}_1-\log(2))\ast S(\mathrm{z}_2^b)+\delta_{b=0}\log(2)S(\mathrm{z}_2^a).
\end{align*}
Then, we obtain that
\begin{align}
\hat{\mathcal{P}}(u,v)&=\mathrm{z}_1\ast \mathcal{S}(v)+\log(2)(\mathcal{S}(u)-\mathcal{S}(v))+\frac{\mathcal{S}(v)}{2}\ast\left[\mathcal{A}\left(\frac{u+v}{2}\right)+\mathcal{A}\left(\frac{u-v}{2}\right)\right]\notag\\
&\qquad+\frac{\mathcal{S}(u)}{2}\ast\left[\mathcal{B}\left(\frac{u+v}{2}\right)+\mathcal{B}\left(\frac{u-v}{2}\right)\right],\label{Eq:hat-P}\\
\hat{\mathcal{P}}^\star(u,v)&=\mathrm{z}_1\ast \mathcal{S}^\star(v)+\log(2)(\mathcal{S}^\star(u)-\mathcal{S}^\star(v))+\frac{\mathcal{S}^\star(v)}{2}\ast\left[\mathcal{A}\left(\frac{u+v}{2}\right)+\mathcal{A}\left(\frac{u-v}{2}\right)\right]\notag\\
&\qquad+\frac{\mathcal{S}^\star(u)}{2}\ast\left[\mathcal{B}\left(\frac{u+v}{2}\right)+\mathcal{B}\left(\frac{u-v}{2}\right)\right].\label{Eq:hat-P-star}
\end{align}
We give a proof of \eqref{Eq:hat-P} and one can prove \eqref{Eq:hat-P-star} similarly. We have $\hat{\mathcal{P}}(u,v)=\hat{\mathcal{P}}_1(u,v)+\hat{\mathcal{P}}_2(u,v)$ with
\begin{align*}
\hat{\mathcal{P}}_1(u,v)&=\sum\limits_{a,b=0}^\infty(-1)^{a+b}\left[\delta_{a=0}(\mathrm{z}_1-\log (2))\ast \mathrm{z}_2^b+\delta_{b=0}\log(2)\mathrm{z}_2^a\right]u^{2a}v^{2b},\\
\hat{\mathcal{P}}_2(u,v)&=\sum\limits_{a,b=0}^\infty\sum\limits_{r=1}^{a+b}(-1)^{a+b+r}2^{-2r}\left[\binom{2r}{2a}+(1-2^{-2r})\binom{2r}{2b}\right]\mathrm{z}_{2r+1}\ast \mathrm{z}_2^{a+b-r}u^{2a}v^{2b}.
\end{align*}
Obviously,
\begin{align*}
\hat{\mathcal{P}}_1(u,v)&=\sum\limits_{b=0}^\infty(-1)^b(\mathrm{z}_1-\log(2))\ast \mathrm{z}_2^bv^{2b}+\sum\limits_{a=0}^\infty(-1)^a\log(2)\mathrm{z}_2^au^{2a}\\
&=\mathrm{z}_1\ast \mathcal{S}(v)-\log(2)\mathcal{S}(v)+\log(2)\mathcal{S}(u).
\end{align*}
Replacing $a+b-r$ by $s$, we get
\begin{align*}
\hat{\mathcal{P}}_2(u,v)&=\sum\limits_{r\geq 1,s\geq0}\sum\limits_{a=0}^{r}(-1)^s2^{-2r}\binom{2r}{2a}\mathrm{z}_{2r+1}\ast \mathrm{z}_2^{s}u^{2a}v^{2s+2r-2a}\\
&\qquad+\sum\limits_{r\geq 1,s\geq0}\sum\limits_{b=0}^{r}(-1)^s2^{-2r}(1-2^{-2r})\binom{2r}{2b}\mathrm{z}_{2r+1}\ast \mathrm{z}_2^{s}u^{2s+2r-2b}v^{2b}\\
&=\sum\limits_{s=0}^\infty(-1)^s\mathrm{z}_2^sv^{2s}\ast\sum\limits_{r=1}^\infty2^{-2r}\mathrm{z}_{2r+1}\sum\limits_{a=0}^{r}\binom{2r}{2a}u^{2a}v^{2r-2a}\\
&\qquad+\sum\limits_{s=0}^\infty(-1)^s\mathrm{z}_2^su^{2s}\ast\sum\limits_{r=1}^\infty2^{-2r}(1-2^{-2r})\mathrm{z}_{2r+1}\sum\limits_{b=0}^{r}\binom{2r}{2b}u^{2r-2b}v^{2b}.
\end{align*}
Using the identity
\begin{align*}
2\sum\limits_{n=0}^r\binom{2r}{2n}u^{2n}v^{2r-2n}=(u+v)^{2r}+(u-v)^{2r},
\end{align*}
we find that
\begin{align*}
\hat{\mathcal{P}}_2(u,v)=\frac{\mathcal{S}(v)}{2}\ast\left[\mathcal{A}\left(\frac{u+v}{2}\right)+\mathcal{A}\left(\frac{u-v}{2}\right)\right]
+\frac{\mathcal{S}(u)}{2}\ast\left[\mathcal{B}\left(\frac{u+v}{2}\right)+\mathcal{B}\left(\frac{u-v}{2}\right)\right],
\end{align*}
which together with the formula of $\hat{\mathcal{P}}_1(u,v)$ deduces \eqref{Eq:hat-P}. And finally, we obtain
\begin{align}\label{Eq:hat-P vs hat-P-star}
\hat{\mathcal{P}}^\star(u,v)=\hat{\mathcal{P}}(v,u)\ast \mathcal{S}^\star(u)\ast \mathcal{S}^\star(v).
\end{align}
Hence, \eqref{Eq:P vs P-star} and \eqref{Eq:hat-P vs hat-P-star} imply the equivalence of the evaluation formulas \eqref{Eq:t_star one-two} and \eqref{Eq:t-one-two}.
\qed

\end{document}